\documentclass[12pt,reqno]{amsart}

\usepackage{amsmath, amsfonts, amsthm, amssymb, graphicx}
\usepackage{float,epsf,subfigure}

\theoremstyle{plain}

\theoremstyle{definition}

\theoremstyle{remark}

\theoremstyle{example}

\numberwithin{equation}{section}

\renewcommand{\r}{\rho}
\newcommand{\g}{\gamma}

\renewcommand{\O}{\Omega}
\newcommand{\G}{\Gamma}
\renewcommand{\k}{\kappa}

\renewcommand{\d}{\partial}

\newcommand{\R}{{\mathbb R}}

\def\be{\begin{equation}}
\def\ee{\end{equation}}
\def\bes{\begin{equation*}}
\def\ees{\end{equation*}}
\def\bc{\begin{cases}}
\def\ec{\end{cases}}

\begin{document}

\title[Fluid Dynamic Formulation for Isometric Embedding]
{A Fluid Dynamic Formulation of the Isometric Embedding Problem in Differential Geometry}

\author{Gui-Qiang Chen}
\address{Department of Mathematics, Northwestern University,
         Evanston, IL 60208, USA.}
\email{gqchen@math.northwestern.edu}

\author{Marshall Slemrod}
\address{Department of Mathematics, University of Wisconsin,
Madison, WI 53706, USA.}
\email{slemrod@math.wisc.edu}

\author{Dehua Wang}
\address{Department of Mathematics, University of Pittsburgh,
                           Pittsburgh, PA 15260, USA.}
\email{dwang@math.pitt.edu}

\date{August 5, 2008} %\today}

\dedicatory{Dedicated to Walter Strauss on the occasion of his 70th birthday}


\begin{abstract}
The isometric embedding problem is a fundamental problem in
differential geometry. A  longstanding problem is considered in this
paper to characterize intrinsic metrics on a two-dimensional
Riemannian manifold which can be realized as isometric immersions
into the three-dimensional Euclidean space. A remarkable connection
between gas dynamics and differential geometry is discussed. It is
shown how the fluid dynamics can be used to formulate a geometry
problem. The equations of gas dynamics are first reviewed. Then the
formulation  using the fluid dynamic variables in conservation laws
of gas dynamics is presented for the isometric embedding problem in
differential geometry.
\end{abstract}
\maketitle

%%%%%%%%%%%%%%%%%%%

\section{Introduction}

We are concerned with isometric embeddings or immersions (i.e.,
realizations) of two-dimensional Riemannian manifolds in the
Euclidean space $\R^3$. A classical question in differential
geometry is whether one can isometrically embed a Riemannian
manifold $({\mathcal M}^n,g)$ into $\R^N$ for $N$ large enough. Nash
\cite{Nash56} indicated that any smooth compact manifold $({\mathcal
M}^n, g)$ can always be isometrically embedded into ${\mathbb{R}}^{3
s_n+4n}$ for $s_n = n(n-1)/2$. This important paper lays  the
foundation for the development of geometric analysis in the second
half of the 20th century. Gromov \cite{Gromov86} proved that one can
embed any $({\mathcal M}^n, g)$ even into
${{\mathbb{R}}}^{s_n+2n+3}$. Then a further natural question is to
find the smallest dimension $N(n)$ for the Riemannian manifold
$({\mathcal M}^n, g)$ to be isometrically embeddable in
${\mathbb{R}}^{N(n)}$. In particular, a fundamental, longstanding
open problem is to characterize intrinsic metrics on a
two-dimensional Riemannian manifold ${\mathcal M}^2$ which can be
realized as isometric immersions into $\R^3$ (cf.
\cite{HanHong,PS,Roz,Yau00} and the references cited therein).
Important results have been achieved for the embedding of surfaces
with positive Gauss curvature which can be formulated as an elliptic
boundary value problem (cf. \cite{HanHong}). For the case of
surfaces of negative Gauss curvature where the underlying partial
differential equations are hyperbolic, the complimentary problem
would be an initial or initial-boundary value problem. When the
Gauss curvature changes sign, the problem then becomes an
initial-boundary value problem of mixed elliptic-hyperbolic type.
Hong in \cite{Hong} first proved that complete negatively curved
surfaces can be isometrically immersed in $\R^3$ if the Gauss
curvature decays at certain rate in the time-like direction. In
fact, a crucial lemma in Hong \cite{Hong} (also see Lemma 10.2.9 in
\cite{HanHong}) shows that, for such a decay rate of the negative
Gauss curvature, there exists a unique global smooth, small solution
forward in time for prescribed smooth, small initial data. We are
interested in solving the corresponding problem for a class of
{large} non-smooth initial data.

In Chen-Slemrod-Wang \cite{CSW2}, we have introduced a general
approach, which combines a fluid dynamic formulation of balance laws
with a compensated compactness framework, to deal with the isometric
immersion problem in $\R^3$ (even when the Gauss curvature changes
sign).  In Chen-Slemrod-Wang \cite{CSW1}, we have developed a
vanishing viscosity method to establish the existence of a weak
entropy solution to the transonic flow in gas dynamics past an
obstacle such as an airfoil, via the method of compensated
compactness (\cite{Mu2, Ta1}). We have found in \cite{CSW2} that the
idea of \cite{CSW1} for gas dynamics is useful for solving the
isometric embedding problem in differential geometry. In particular,
in \cite{CSW2}, we have formulated the isometric immersion problem
for two-dimensional Riemannian manifolds in $\R^3$ via solvability
of the Gauss-Codazzi system, and have introduced a fluid dynamic
formulation of balance laws for the Gauss-Codazzi system. Then we
have formed a compensated compactness framework and present one of
our main observations that this framework is a natural formulation
to ensure the weak continuity of the Gauss-Codazzi system for
approximate solutions, which yields the isometric realization of
two-dimensional surfaces in $\R^3$. As a first application of this
approach,  we have focused on the isometric immersion problem of
two-dimensional Riemannian manifolds with strictly negative Gauss
curvature. Since the local existence of smooth solutions follows
from the standard hyperbolic theory, we are concerned with the
global existence of solutions of the initial value problem with
large initial data. The metrics $(g_{ij})$ we study have special
structures and forms usually associated with  the catenoid of
revolution and the helicoid. For these cases, while Hong's theorem
\cite{Hong} applies to obtain the existence of a solution for small
smooth initial data, our result yields a large-data existence
theorem for a $C^{1,1}$ isometric immersion. To achieve this, we
have introduced a vanishing viscosity method depending on the
features of the initial value problem for isometric immersions and
have presented a technique to make the apriori estimates including
the $L^\infty$ control and $H^{-1}$--compactness for the viscous
approximate solutions. This yields the weak convergence of the
vanishing viscosity approximate solutions and the weak continuity of
the Gauss-Codazzi system for the approximate solutions, hence the
existence of a $C^{1,1}$--isometric immersion of the manifold into
$\R^3$ with prescribed initial conditions.

From Chen-Slemrod-Wang \cite{CSW1, CSW2}, we have seen a remarkable
connection between the two distinct areas of gas dynamics and
differential geometry. Here we present such a connection and show
how the fluid dynamics can be used to formulate a geometry problem.
Thus, we will present first the equations in  Chen-Slemrod-Wang
\cite{CSW1} for  the transonic flow problem in gas dynamics, and
then the formulation using the fluid dynamic variables in
conservation laws of gas dynamics for the isometric embedding
problem in differential geometry.

%%%%%%%%%%%%%%%%%
\section{Equations of Gas Dynamics}

In two space dimensions with variables $(x,y)$, the steady transonic flow
of isentropic case is governed by the following steady Euler equations
on conservations of mass and momentum in gas dynamics:
\begin{equation} \label{euler1}
\begin{cases}
(\rho u)_x + (\rho v)_y = 0,   \\
(\rho u^2 + p)_x + (\rho uv)_y = 0,\\
(\rho u v)_x + (\rho v^2 + p)_y = 0,
\end{cases}
\end{equation}
where $\rho$ is the density, $(u,v)$ is the velocity, and
 $p = \frac{\rho^\gamma}{\gamma} \; (\gamma \geq 1)$ is the pressure.
If we assume that the flow is irrotational, then system
\eqref{euler1} can be reduced to the following two equations of
irrotationality and conservation of mass:
\begin{equation} \label{euler2}
\begin{cases}
v_x - u_y = 0,  \\
(\rho u)_x + (\rho v)_y = 0,
\end{cases}
\end{equation}
and, by scaling, the density $\rho$ is determined by Bernoulli's
law:
\begin{equation} \label{bernoulli}
 \rho = \Big( 1 - \frac{\gamma -1} 2  q^2 \Big)^{1\over \gamma - 1},
\end{equation}
where $q$ is the flow speed defined by $q^2 =u^2 + v^2$.
The sound speed $c$ is defined as
\begin{equation} \label{sp}
c^2=p'(\r)=1-\frac{\g-1}2q^2.
\end{equation}
At the cavitation point $\r=0$,
 $$q=q_{cav}:= \sqrt{\frac2{\g-1}}.$$
 At the stagnation point $q=0$, the density reaches its maximum
 $\r=1$. Bernoulli's law \eqref{bernoulli} is valid for $0\le q\le
 q_{cav}$.  At the sonic point
$q=c$, \eqref{sp} implies $q^2=\frac2{\g+1}$. Define the critical
speed $q_{cr}$ as
$$q_{cr}:=\sqrt{\frac2{\g+1}}.$$
We rewrite Bernoulli's law \eqref{bernoulli} in the form
\begin{equation} \label{ber9}
q^2-q_{cr}^2=\frac2{\g+1}\left(q^2-c^2\right).
\end{equation}
 Thus the flow is subsonic when $q<q_{cr}$, sonic when $q=q_{cr}$,
 and supersonic when $q>q_{cr}$.
For the isothermal flow ($\g=1$), $p={c}^2\r$ where ${c}>0$ is the
constant sound speed,  the density $\r$ is given by Bernoulli's
law:
\begin{equation} \label{ber2}
\r=\r_0\exp\big(-\frac{u^2+v^2}{2{c}^2}\big)
\end{equation}
for some constant $\r_0>0$, and  $q_{cr}={c}$.

%%%%%%%%%%%%%%%%

\section {Isometric Embedding in Differential Geometry}

In this section, we discuss the isometric embedding problem in
differential geometry in $\R^3$ and its formulation of fluid
dynamics.

We first give the Gauss-Codazzi system of isometric embedding in
$\R^3$. Let $g_{ij}, i, j=1,2,$ be the given metric of a
two-dimensional Riemannian manifold $\mathcal{M}$ parameterized on
an open set $\O\subset\R^2$. The first fundamental form $I$ for
$\mathcal{M}$ on $\Omega$ is
\begin{equation*} %\label{2.3}
I:=g_{11}(dx)^2+2 g_{12}dxdy +g_{22}(dy)^2,
\end{equation*}
and {the isometric embedding problem} is to {seek a map
$\mathbf{r}: \Omega\to \R^3$ such that
$
d\mathbf{r}\cdot d\mathbf{r}=I,
$
that is,
\begin{equation*} %\label{2.4}
\partial_x{\bf r}\cdot\partial_x{\bf r}=g_{11},\quad
    \partial_x{\bf r}\cdot\partial_y{\bf r}=g_{12},\quad
   \partial_y{\bf r}\cdot\partial_y{\bf r}=g_{22},
\end{equation*}
so that $\{\partial_x{\bf r}, \partial_y{\bf r}\}$ in $\R^3$ are
linearly independent}.
The corresponding second fundamental form is
\begin{equation*}%\label{2.5}
I\!I:=h_{11}(dx)^2+ 2h_{12}dxdy + h_{22}(dy)^2.
\end{equation*}
The fundamental theorem of surface theory (cf.
\cite{doC1992,HanHong}) indicates that { there exists a surface in
$\R^3$ whose first and second fundamental forms are $I$ and $I\!I$
if the coefficients $(g_{ij})$ and $(h_{ij})$ of the two given
quadratic forms $I$ and $I\!I$ with $(g_{ij})>0$ satisfy the
Gauss-Codazzi system}. It is indicated in Mardare \cite{Mardare2}
(Theorem 9; also see \cite{Mardare1}) that this theorem holds even
when $(h_{ij})$ is only in $L^\infty$ for given $(g_{ij})$ in
$C^{1,1}$, for which the immersion surface is $C^{1,1}$. This shows
that, for the realization of a two-dimensional Riemannian manifold
in $\R^3$ with given metric $(g_{ij})>0$, it suffices to solve
$(h_{ij})\in L^\infty$ determined by the Gauss-Codazzi system to
recover ${\bf r}$ a posteriori. The Gauss-Codazzi system (cf.
\cite{doC1992,HanHong}) can be written as
\begin{equation} \label{g1}
\begin{cases}
\d_x{M}-\d_y{L}=\G^{(2)}_{22}L-2\G^{(2)}_{12}M+\G^{(2)}_{11}N, \\
\d_x{N}-\d_y{M}=-\G^{(1)}_{22}L+2\G^{(1)}_{12}M-\G^{(1)}_{11}N,
\end{cases}
\end{equation}
with
\begin{equation}\label{g2}
LN-M^2=\k,
\end{equation}
where
$$
L=\frac{h_{11}}{\sqrt{|g|}}, \quad M=\frac{h_{12}}{\sqrt{|g|}},
\quad N=\frac{h_{22}}{\sqrt{|g|}},\quad
|g|=det(g_{ij})=g_{11}g_{22}-g_{12}^2, $$
 $\kappa(x,y)$ is the Gauss
curvature that is determined by the relation:
$$
\k(x,y)=\frac{R_{1212}}{|g|}, \;
R_{ijkl}=g_{lm}\left(\d_k\G^{(m)}_{ij}-\d_j\G^{(m)}_{ik}
+\G^{(n)}_{ij}\G^{(m)}_{nk}-\G^{(n)}_{ik}\G^{(m)}_{nj}\right),
$$
$R_{ijkl}$ is the curvature tensor and depends on $(g_{ij})$ and its
first and second derivatives, and
$$
\G_{ij}^{(k)}=\frac12g^{kl}\left(\d_j g_{il}+\d_i g_{jl}-\d_l
 g_{ij}\right)
$$
is the Christoffel symbol and depends on the first derivatives of
$(g_{ij})$, where the summation convention is used, $(g^{kl})$
denotes the inverse of $(g_{ij})$, and
$(\partial_1,\partial_2)=(\partial_x, \partial_y)$.
Therefore, given a positive definite metric $(g_{ij})\in C^{1,1}$,
the Gauss-Codazzi system gives us three equations for the three
unknowns $(L, M, N)$ determining the second fundamental form $I\!I$.
Note that, although $(g_{ij})$ is positive definite, $R_{1212}$ may
change sign and so does the Gauss curvature $\k$. Thus, the
Gauss-Codazzi system \eqref{g1}--\eqref{g2} generically is of mixed
hyperbolic-elliptic type, as in transonic flow (cf. \cite{CSW1}). In
Chen-Slemrod-Wang \cite{CSW2}, we have introduced a general approach
to deal with the isometric immersion problem involving nonlinear
partial differential equations of mixed hyperbolic-elliptic type by
combining a fluid dynamic formulation of balance laws  with a
compensated compactness framework. As an example of direct
applications of this approach, we have shown how this approach can
be applied to establish an isometric immersion of a two-dimensional
Riemannian manifold with negative Gauss curvature in $\R^3$.

We now describe the fluid dynamic formulation of the Gauss-Codazzi
system \eqref{g1}--\eqref{g2} in detail. Although, from the
viewpoint of geometry, the constraint condition \eqref{g2} is a
Monge-Amp\`{e}re equation and the equations in \eqref{g1} are
integrability relations, we can put the problem into a fluid dynamic
formulation so that the isometric immersion problem may be solved
via the approaches for transonic flows of fluid dynamics in
Chen-Slemrod-Wang \cite{CSW1}. To do this, we set
$$
L=\r v^2+p,  \quad M=-\r uv, \quad  N=\r u^2+p,
$$
and  set $q^2=u^2+v^2$ as usual. Then the equations in \eqref{g1}
become the familiar balance laws of momentum:
\begin{equation} \label{g3}
\begin{cases}
\d_x(\r uv)+\d_y(\r v^2+p)
 =-(\r v^2+p)\G^{(2)}_{22}-2\r uv\G^{(2)}_{12}-(\r u^2+p)\G^{(2)}_{11}, \\
\d_x(\r u^2+p)+\d_y(\r uv)
 =-(\r v^2+p)\G^{(1)}_{22}-2\r uv\G^{(1)}_{12}-(\r u^2+p)\G^{(1)}_{11},
\end{cases}
\end{equation}
and the Monge-Amp\`{e}re constraint \eqref{g2} becomes
\begin{equation}\label{g4}
\r p q^2+p^2=\k.
\end{equation}
We choose pressure $p$ as for the Chaplygin-type gas:
\begin{equation*}\label{g5}
p=-\frac{1}{\r}.
\end{equation*}
Then, from \eqref{g4},  we have the ``Bernoulli" relation:
\begin{equation}\label{g6}
\r=\frac{1}{\sqrt{q^2+\k}}.
\end{equation}
This yields
\begin{equation}\label{g7}
p=-\sqrt{q^2+\k},
\end{equation}
and the formulas for $u^2$ and $v^2$:
$$
u^2=p(p-M), \qquad v^2=p(p-L), \qquad M^2=(N-p)(L-p).
$$
The last relation for $M^2$ gives the relation for $p$ in terms of
$(L,M,N)$, and then the first two give the relations for $(u, v)$ in
terms of $(L,M,N)$.

We rewrite \eqref{g3} as
\begin{equation} \label{g8}
\begin{cases}
\d_x(\r uv)+\d_y(\r v^2+p) =R_1, \\
\d_x(\r u^2+p)+\d_y(\r uv) =R_2,
\end{cases}
\end{equation}
where $R_1$ and $R_2$ denote the right-hand sides of \eqref{g3}.
Then we can write down our ``rotationality-continuity equations" as
\begin{equation} \label{g11}
\begin{cases}
\d_x v- \d_y u =\frac1{\r q^2}\Big(u\big(\frac12\r \d_y\k+R_1\big)-
v\big(\frac12\r \d_x\k+R_2\big)\Big),  \\
\d_x(\r u)+ \d_y(\r v)= \frac12\frac{\r u}{q^2} \d_x\k
+\frac12\frac{\r v}{q^2} \d_y\k
 +\frac{v}{q^2}R_1+\frac{u}{q^2}R_2.
\end{cases}
\end{equation}
In summary, the Gauss-Codazzi system \eqref{g1}--\eqref{g2}, the
momentum equations \eqref{g3}--\eqref{g7},
and the rotationality-continuity equations \eqref{g6} and
\eqref{g11} are all formally equivalent.
However, for weak solutions, we know from our experience with gas
dynamics that this equivalence breaks down. In
Chen-Dafermos-Slemrod-Wang \cite{CDSW}, the decision has been made
(as is standard in gas dynamics) to solve the
rotationality-continuity equations and view the momentum equations
as ``entropy" equalities which may become inequalities for weak
solutions. In geometry, this situation is just the reverse. It is
the Gauss-Codazzi system that must be solved exactly, and hence the
rotationality-continuity equations will become ``entropy"
inequalities for weak solutions.

We define the ``sound" speed as:
\begin{equation}\label{g18}
c^2=p'(\r)=\frac1{\r^2},
\end{equation}
then from our ``Bernoulli" relation  \eqref{g6}, we see
\begin{equation}\label{g20}
 c^2=q^2+\k.
\end{equation}
Hence, under this formulation,

(i) when $\k>0$, the ``flow" is subsonic, i.e., $q<c$,
       and system \eqref{g3}--\eqref{g4} is elliptic;

(ii) when $\k<0$, the ``flow" is supersonic, i.e., $q>c$,
       and system \eqref{g3}--\eqref{g4} is hyperbolic;

(iii) when $\k=0$, the ``flow" is sonic, i.e., $q=c$,
       and system \eqref{g3}--\eqref{g4} is degenerate.

In general, system \eqref{g3}--\eqref{g4} is of mixed
hyperbolic-elliptic type. Thus, the isometric immersion problem
involves the existence of solutions to nonlinear partial
differential equations of mixed hyperbolic-elliptic type.

In Chen-Slemrod-Wang \cite{CSW2}, we have considered one of the
spatial variables $x$ and $y$ as time-like, have introduced a
vanishing viscosity method via parabolic regularization to obtain
the uniform $L^\infty$ estimate by identifying invariant regions for
the approximate solutions, and have shown that the
$H^{-1}_{loc}$--compactness can be achieved for the viscous
approximate solutions. Then, as in Chen-Slemrod-Wang \cite{CSW1},
the compensated compactness framework yields a weak solution to the
initial value problem of system \eqref{g3}--\eqref{g4} when the
initial data lies in the diamond-shaped invariant region. This
establishes a  $C^{1,1}(\R^2)$ immersion of the Riemannian manifold
into $\R^3$. In particular, our existence result asserts the
existence of a $C^{1,1}$-surface for the associated metric for a
class of non-circular cross-sections prescribed at $x=0$ for
catenoid. Our study in \cite{CSW2} also applies to the helicoid. See
\cite{CSW2} for the details. Possible implication of our approach
may be in existence theorems for equilibrium configurations of a
catenoidal shell as detailed in Vaziri-Mahedevan \cite{VM}. However,
the existence of isometric embeddings/immersions of a general
surface with negative Gauss curvature is still open. When  the Gauss
curvature $\k$ changes sign, the problem becomes transonic and thus
mixed hyperbolic-elliptic type. In this mixed-type problem, only
special local solutions are known to exist for special data
(\cite{CSLin,HanHong}), and the existence of global solutions is a
significantly difficult open problem.

\bigskip\bigskip
{\bf Acknowledgments.} Gui-Qiang Chen's research was supported in part by the National
 Science Foundation under Grants DMS-0807551, DMS-0720925,
 and DMS-0505473.  Marshall Slemrod's research was supported in part by
 the National Science Foundation under Grant DMS-0647554.
 Dehua Wang's research was supported in part by the National Science
Foundation under Grant DMS-0604362, and by the Office of Naval Research
under Grant N00014-07-1-0668.

\end{document}